\documentclass[a4paper,11pt]{article}
\usepackage[T1]{fontenc}
\usepackage{latexsym}
\usepackage{amssymb,amsmath,mathtools,theorem}
\usepackage[a4paper]{geometry}
\newtheorem{thm}{Theorem}[section]
\newtheorem{lem}[thm]{Lemma}
\newtheorem{prop}[thm]{Proposition}

{\theorembodyfont{\rmfamily}  \newtheorem{defn}[thm]{Definition}}
{\theorembodyfont{\rmfamily}  \newtheorem{rem}[thm]{Remark}}
{\theorembodyfont{\rmfamily}  \newtheorem{exam}[thm]{Example}}
\newcommand{\J}{\mathcal J}
\newcommand{\ep}{\varepsilon}
\newcommand{\C}{\mathbb{C}}

\newcommand{\N}{\mathbb{N}}

\newcommand{\f}{\varphi}
\newcommand{\re}{\mathrm{Re}}

\newcommand{\lm}{\lambda}

\newcommand{\T}{\mathcal T}

\newcommand{\Hh}{\mathcal H}

\newcommand{\WAP}{\operatorname{WAP}}
\newcommand{\WM}{\operatorname{WM}}

\newcommand{\Cbr}{C_\mathrm{b,r}}
\newcommand{\Cb}{C_\mathrm{b}}
\newcommand{\mf}{\mathfrak{m}}

\newcommand{\la}{\langle}
\newcommand{\ra}{\rangle}

\newcommand{\Span}{\operatorname{Span}}

\newcommand{\ov}{\overline}
\newcommand{\wh}{\widehat}
\newcommand{\CP}{\mathrm{CP}}
\newcommand{\Tr}{\mathrm{Tr}}

\title{\textbf{Non Kazhdan's groups and a new approximation property for tracial von Neumann algebras}}
\author{Paul Jolissaint}

\begin{document}

\maketitle

\begin{abstract}
We define a new approximation property for tracial von Neumann algebras, called \textit{weakly mixing approximation property} which, for discrete groups and II$_1$ factors, is equivalent to the negation of Kazhdan's property (T).
\end{abstract}

\medskip\noindent
\emph{Mathematics Subject Classification:} 46L10, 22D55\\
\emph{Key words:} Weakly mixing representations, weakly almost periodic functions, tracial von Neumann algebras, Kazhdan's property (T).

\section{Introduction}

For several years, many approximation properties for discrete groups have been translated with success to their associated von Neumann algebras; it was first the case for amenability, then for the Haagerup property, and for Kazhdan's property (T), starting in \cite{CoJo}. More precisely, property (T) for finite von Neumann algebras defined in the latter article is equivalent to Kazhdan's property (T) for groups with finitely many finite conjugacy classes only, see \cite[Theorem A]{JolT}. Next, equivalent definitions of property (T) for finite von Neumann algebras (and for pairs of such algebras) were defined and studied rather extensively: see for instance \cite{AS} and \cite{Popa}.
Thus, we think that it is interesting to have another characterization for its negation, and that is what is presented in our note.

Our approach is based on the notion of coefficients of unitary representations on bimodules coming from completely positive maps on tracial von Neumann algebras; as these coefficients are weakly almost periodic functions (as are all coefficients of unitary representations), this allows to consider weakly mixing representations (see Section 2). Thus, we introduce an approximation property for tracial von Neumann algebras, called \textit{weakly mixing approximation property} (WMAP), which means that, for a given tracial von Neumann algebra $(M,\tau)$, there exists a net of completely positive maps $(\phi_i)$ on $M$ such that $\lim_i \Vert \phi_i(x)-x\Vert_2=0$ for every $x\in M$ and such that the unitary representation $\pi_{i}$ of the unitary group $U_M$ on the $M$-bimodule $\Hh_{\phi_i}$ defined by $\pi_i(u)\xi=u\xi u^*$ is weakly mixing; see Definition \ref{defwmap}. 

Although we had hoped to have discovered a complete characterization of tracial von Neumann algebras which do not have property (T), we realized that it holds only for infinite discrete groups and their associated von Neumann algebras, and for II$_1$ factors. In fact, we give below a rather simple example of a tracial von Neumann algebra that has neither property (T) nor the WMAP.

\medskip
Here is our main result:

\begin{thm}\label{mainthm}
\begin{enumerate}
\item [(1)] If $(M,\tau)$ is a tracial von Neumann algebra that has the WMAP, and if $\tau'$ is another normal, faithful tracial state on $M$, then $(M,\tau')$ has the same property.
\item [(2)] If $G$ is a discrete group, then it does not have Kazhdan's property (T) if and only if its associated von Neumann algebra $L(G)$ has the WMAP.
\item [(3)] If $(M,\tau)$ is a tracial von Neumann algebra that has the WMAP, then it does not have property (T).
\item [(4)] If moreover $(M,\tau)$ is $\mathrm{II}_1$ factor, then it has the WMAP if and only if it does not have property (T).
\end{enumerate}
\end{thm}
As will be reminded in Section 3, every finite, atomic von Neumann algebra with separable predual has property (T), and thus it does not have the WMAP.

\medskip
Here is the promised example of a tracial, diffuse von Neumann algebra that has neither property (T) nor the WMAP:

\begin{exam}
Let $A$ be a II$_1$ factor with property (T) and let $B$ be a II$_1$ factor which does not have property (T). Set $M=A\oplus B$. Then, as property (T) is stable by reduction by nonzero projections of $M$ (\cite[Proposition 4.7]{Popa}), $M$ does not have property (T). As is proven in Proposition \ref{HAP}, the WMAP is also preserved by reduction by nonzero central projections. Hence by Theorem \ref{mainthm}(3), $M$ cannot have the WMAP either.
\end{exam}

\medskip
The rest of the article is organized as follows: Section 2 contains known material on weakly mixing representations that will be needed in Section 3, the latter being devoted to the definition and study of WMAP, including the proof of the main results. Our undefined notations stick to those of the monograph \cite{AS}, except that our scalar products are linear in the first variable.

\medskip\noindent
\textit{Acknowledgements.} I am very grateful to the referee for his/her very pertinent remarks and valuable suggestions, and in particular for having detected a gap in the proof of Theorem \ref{final} in a preliminary version.

\section{Unitary representations of topological groups, non Kazhdan groups and weak mixing}

We gather in this section known definitions and results that will be useful in Section 3.
\medskip

Let $G$ be a topological group, which is assumed to be Hausdorff throughout the article, and denote by $e$ the identity element of $G$. Let $\Cb(G)$ (resp. $\Cbr(G)$) be the commutative $C^*$-algebra of all bounded, continuous (resp. uniformly right continuous)
functions on $G$ equipped with the uniform norm $\Vert f\Vert_\infty\coloneqq \sup_{s\in G}|f(s)|$. As in \cite{Jol1}, for $g\in G$ and $f:G\rightarrow \C$, we denote by $g\cdot f:G\rightarrow\C$ (resp. $f\cdot g$) the left (resp. right) translate of $f$ by $g$, i.e.
\[
(g\cdot f)(s)=f(g^{-1}s)\quad\textrm{and}\quad
(f\cdot g)(s)=f(sg)
\]
for all $f:G\rightarrow\C$ and $g,s\in G$, the corresponding left (resp. right) orbit being denoted by $Gf$ (resp. $fG$).

Thus, a function $f\in\Cb(G)$ belongs to $\Cbr(G)$ if and only if $\Vert g\cdot f-f\Vert_\infty\to 0$ as $g\to e$.

A function $f\in \Cbr(G)$ is \textit{weakly almost periodic} if its orbit $Gf$ is weakly relatively compact in $\Cbr(G)$. It turns out that if $f\in \Cbr(G)$ is weakly almost periodic, then its right orbit $fG$ is also weakly relatively compact (\cite[Corollary 1.12]{Bur}).

We denote by $\WAP(G)$ the sub-$C^*$-algebra of weakly almost periodic functions on $G$. (We observe that the right uniform continuity of almost periodic functions is used to prove that $\WAP(G)$ is indeed a $C^*$-algebra: see the proof of \cite[Proposition 5.4]{Jol1}.) 

The main feature of $\WAP(G)$ is that it admits a unique bi-invariant mean (see for instance \cite[Chapter 1]{Bur} or  \cite[Theorem 5.5]{Jol1}), more precisely:

\begin{thm}\label{moy}
Let $G$ be a topological group. There exists a unique linear functional $\mf:\WAP(G)\rightarrow\C$ with the following properties:
\begin{enumerate}
\item [(i)] $\mf(f)\geq 0$ for every $f\geq 0$;
\item [(ii)] $\mf(1)=1$;
\item [(iii)] $\mf(g\cdot f)=\mf(f\cdot g)=\mf(f)$ for all $f\in\WAP(G)$ and 
$g\in G$;
\item [(iv)] for every $f\in\WAP(G)$ and every $\ep>0$, there exists a convex combination $\displaystyle{\psi:=\sum_{j=1}^m t_j{g_j}\cdot f}$ (with $g_j\in G$ and $t_j\geqslant 0$, $\sum\limits_j t_j=1$) such that 
$
\Vert\psi-\mf(f)\Vert_\infty<\ep,
$
and there exists a convex combination $\displaystyle{\f:=\sum_i s_if\cdot h_i}$ (with $h_i\in G$ and $s_i\geqslant 0$, $\displaystyle{\sum_i s_i=1}$) such that
$
\Vert \f-\mf(f)\Vert_\infty<\ep.
$
\end{enumerate}
\end{thm}

The following result is used in the next section, and as we do not know if it is already known, we provide a proof.

\begin{prop}\label{sbgroup}
Let $G$ be a topological group and $H$ a closed subgroup of $G$. Let us denote by $\mf_G$ (resp. $\mf_H$) the bi-invariant mean on $\WAP(G)$ (resp. $\WAP(H)$). Then $\mf_H(f|_H)=\mf_G(f)$ for every $f\in \WAP(G)$.
\end{prop}
\textsc{Proof.} Let us fix $f\in \WAP(G)$; replacing it by $f-\mf_G(f)$, we assume that $\mf_G(f)=0$. Referring to the proof of \cite[Theorem 5.5]{Jol1}, let us set the following notations:
\begin{itemize}
\item $Q_{\mathrm l}(f,G)$ (resp. $Q_{\mathrm r}(f,G)$) denotes the norm-closed convex hull of the left orbit $Gf$ (resp. right orbit $fG$) in $\Cbr(G)$.
\item Similarly, $Q_{\mathrm l}(f|_H,H)$ and $Q_{\mathrm r}(f|_H,H)$ denote the corresponding sets in $\Cbr(H)$.
\item Finally, $K_{\mathrm l}(f)$ (resp. $K_{\mathrm r}(f)$) denotes the norm-closed convex hull of $Hf$ (resp. $fH$) in $\Cbr(G)$.
\end{itemize}
Then obviously, $K_{\mathrm l}(f)\subset Q_{\mathrm l}(f,G)$ and $K_{\mathrm r}(f)\subset Q_{\mathrm r}(f,G)$, so that $K_{\mathrm l}(f)$ and $K_{\mathrm r}(f)$ are weakly compact. As the left and right actions of $H$ on $K_{\mathrm l}(f)$ and $K_{\mathrm r}(f)$ respectively are obviously distal, by Ryll-Nardzewski Theorem, $K_{\mathrm l}(f)$ and $K_{\mathrm r}(f)$ both contain constant functions, and, by a standard convexity argument, these constant functions are unique and equal, which is denoted by $c(f)$. As $K_{\mathrm l}(f)\cap K_{\mathrm r}(f)\subset Q_{\mathrm l}(f,G)\cap Q_{\mathrm r}(f,G)$, one has necessarily $c(f)=\mf_G(f)=0$.

Let $\ep>0$; there exists $t_1,\ldots,t_m>0$, $\sum_j t_j=1$, and $h_1,\ldots,h_m\in H$ such that
\[
\sup_{g\in G}\left|\sum_j t_j f(h_j^{-1}g)\right|<\ep.
\]
A fortiori,
\[
\sup_{h\in H}\left|\sum_j t_j f(h_j^{-1}h)\right|<\ep,
\]
which shows that $0\in Q_{\mathrm l}(f|_H,H)\cap Q_{\mathrm r}(f|_H,H)$ hence that $\mf_H(f|_H)=0$.
\hfill $\square$

\medskip
A \textit{unitary representation} of $G$ (often called a \textit{representation} hereafter for short) is a pair $(\pi,\Hh)$ where $\pi:G\rightarrow U(\Hh)$ is a continuous homomorphism from $G$ to the unitary group of $\Hh$ equipped with the strong operator topology. This means that, for every $\xi\in\Hh$, the map $g\mapsto \pi(g)\xi$ is continuous.

Given such a representation and vectors $\xi,\eta\in \Hh$, the associated \textit{coefficient function} $\xi*_\pi\bar{\eta}$ is defined on $G$ by
\[
\xi*_\pi\bar{\eta}(g)\coloneqq \la \pi(g)\xi|\eta\ra \quad (g\in G).
\]
One has for all $g,h\in G$:
\[
g\cdot(\xi*_\pi\bar\eta)\cdot h=(\pi(h)\xi)*_\pi(\overline{\pi(g)\eta})
\]
and it is well known and easy to verify that every coefficient function belongs to $\Cbr(G)$, and even to $\WAP(G)$; the latter assertion follows from the weak compactness of the orbit $\{\pi(g)\eta\colon g\in G\}$ in $\Hh$ and from the continuity of the map $\pi(g)\eta\mapsto\xi*_\pi\ov{\pi(g)\eta}=g\cdot(\xi*_\pi\bar\eta)$ with respect to the weak topologies (see for instance \cite[Theorem 3.1]{Bur}).

If $(\pi,\Hh)$ is a unitary representation of $G$, then the function $\xi*_\pi\bar\xi$ is of positive type for every $\xi\in\Hh$, and conversely, to every function of positive type $\f$, the GNS construction associates a unique triple (up to unitary isomorphism) $(\pi_\f,\Hh_\f,\xi_\f)$ where $(\pi_\f,\Hh_\f)$ is a unitary representation of $G$, $\xi_\f\in \Hh_\f$ is a cyclic vector, i.e. the orbit $\pi_\f(G)\xi_\f$ is a total subset of $\Hh_\f$, and
$\f=\xi_\f*_{\pi_\f}\ov{\xi_\f}$. See for instance \cite[Theorem C.4.10]{BHV}.

\medskip
For later use, given a representation $(\pi,\Hh)$ and $\xi,\eta\in \Hh$, let us recall the following useful formula for $\mf(\xi*_\pi\bar\eta)$:
denote by $\Hh^\pi$ the subspace of invariant vectors under the action of $\pi(G)$, and by $P_\pi$ the orthogonal projection from $\Hh$ onto $\Hh^\pi$. Then one has
\begin{equation}\label{proj}
    \mf(\xi*_\pi\bar\eta)=\la P_\pi\xi|\eta\ra=\la \xi|P_\pi\eta\ra.
\end{equation}
It is a consequence of Theorem \ref{moy}(iv) and of the fact that $P_\pi\zeta$ is the element of minimal norm in the norm-closed convex hull of $\pi(G)\zeta$ for every $\zeta\in\Hh$.

\begin{defn}\label{def2.3}
Let $G$ be a topological group.
\begin{enumerate}
\item [(1)] (\cite[Definition 1.1]{BR}) A representation $(\pi,\Hh)$ of $G$ is \textit{weakly mixing} if for all $\xi,\eta\in \Hh$, one has $\mf(|\xi*_\pi\bar\eta|)=0$.
\item [(2)] Let $(\pi,\Hh)$ be a unitary representation of $G$. Then $\pi$ \textit{almost has invariant vectors} if, for every compact subset $Q\not=\emptyset$ of $G$ and every $\ep>0$, there exists a unit vector $\xi\in\Hh$ such that
\[
\sup_{g\in Q}\Vert \pi(g)\xi-\xi\Vert<\ep.
\]
We denote that property by $1_G\prec\pi$, as it corresponds to the weak containment of the trivial representation in the case of locally compact groups.
\item [(3)] (\cite[Definition 1.1(2)]{Jol2}) The topological group $G$ has \textit{Kazhdan's property (T)} if, for every unitary representation $(\pi,\Hh)$ of $G$ such that $1_G\prec \pi$, there exists a nonzero vector $\xi\in\Hh$ such that $\pi(g)\xi=\xi$ for every $g\in G$.
\end{enumerate}
\end{defn}
As one more application of Theorem \ref{moy}(iv), it is easy to see that, for $\xi,\eta\in \Hh$, one has $\mf(|\xi*_\pi\bar\eta|)=0$ if and only if $\mf(|\xi*_\pi\bar\eta|^2)=0$, and it is
obvious that, if a representation contains a nonzero, invariant vector $\xi$, then the associated coefficient is a nonzero constant, hence the representation cannot be weakly mixing. This implies that if $G$ has a weakly mixing representation that almost has invariant vectors, then $G$ does not have property (T).

Conversely, if $G$ is moreover locally compact and $\sigma$-compact, then 
the converse is true, as the following result shows (see \cite[Theorem 1.2]{Jol2} where the result is implicit). In case $G$ is moreover separable, see also \cite[Theorem 1]{BV}. In fact, it rests essentially on the fact that a locally compact, $\sigma$-compact group $G$ does not have property (T) if and only if it admits an unbounded real-valued function conditionally of negative type: see \cite[Lemma 4.4]{JolT} and the proof of \cite[Theorem 1.2]{Jol2}.

\begin{thm}
Suppose that $G$ is a locally compact, $\sigma$-compact group. Then it admits a weakly mixing representation such that $1_G\prec\pi$ if and only if it does not have Kazhdan's property (T).
\end{thm}
Given $f\in\WAP(G)$, it can be somewhat complicated to verify that $\mf(|f|)=0$. The following criterium for coefficients of representations will be used here. It is essentially contained in \cite[Theorems 1.3 and 1.9]{BR}, which are stated and proved there for locally compact, $\sigma$-compact groups, but their proofs hold for arbitrary topological groups. We provide here an essentially self-contained proof for the sake of completeness and for the reader's convenience.

\begin{prop}\label{Prop2.4}
Let $G$ be a topological group, $(\pi,\Hh)$ a unitary representation of $G$ and $\T$ any total subset of $\Hh$. Then the following conditions are equivalent:
\begin{enumerate}
\item [(a)] [respectively (a')] the representation $\pi$ is weakly mixing [respectively, for all $\xi,\eta\in \T$, one has $\mf(|\xi*_\pi\bar\eta|)=0$];
\item [(b)] [respectively (b')] for all vectors $\xi,\eta\in \Hh$ [respectively for all $\xi,\eta\in \T$], for every finite set $S\subset G$ and for every $\ep>0$, there exists $g\in G$ such that 
\[
\max_{s\in S}|\la\pi(gs)\xi|\eta\ra|<\ep;
\]
\item [(c)] the representation $\pi$ does not contain any nonzero finite-dimensional subrepresentation.
\end{enumerate}
\end{prop}
\textsc{Proof.} First, it follows from property (iv) of Theorem \ref{moy} that, if $f\in\WAP(G)$ is a real-valued function, then 
\begin{equation}\label{eq1}
 \inf_{g\in G}f(g)\leq \mf(f)\leq\sup_{g\in G}f(g).
\end{equation}
Equivalences between (a) and (a') on the one hand, and of (b) and (b') on the other hand, follow from the positivity of the mean $\mf$, from standard approximation arguments and from the following general inequality valid for all vectors $\xi,\xi',\eta,\eta'\in \Hh$:
\[
\Vert \xi*_\pi,\bar\eta-\xi'*_\pi\bar\eta'\Vert_\infty\leq \Vert \xi-\xi'\Vert\Vert\eta\Vert+\Vert\xi'\Vert\Vert\eta-\eta'\Vert.
\]
$(a)\ \Rightarrow\ (b)$. Assume that $\pi$ is weakly mixing and let $\xi,\eta,S$ and $\ep>0$ be as in (b). By invariance of $\mf$, the equality $\mf(|\xi*_\pi\bar\eta|)=0$ implies that 
\[
\mf\Big(\sum_{s\in S} |\xi*_\pi\bar\eta|\cdot s\Big)=
\mf\Big(\sum_{s\in S}|(\pi(s)\xi)*_\pi\bar\eta|\Big)=0
\]
and thus, by the first inequality in (\ref{eq1}), 
\[
\inf_{g\in G}\Big(\sum_{s\in S}|\la\pi(gs)\xi|\eta\ra|\Big)=0
\]
which implies in turn the existence of some $g\in G$ such that
\[
\sum_{s\in S}|\la\pi(gs)\xi|\eta\ra|<\ep.
\]
$(b)\ \Rightarrow\ (a)$. Let $\xi,\eta\in\Hh$ and $\ep>0$. Set $f=|\xi*_\pi\bar\eta|$. Again by Theorem \ref{moy}(iv), there exists $g_1,\ldots,g_n\in G$ and $t_1,\ldots,t_n\geq 0$, $\sum_j t_j=1$, such that 
\[
\sup_{g\in G}\left|\sum_j t_jf(g_j^{-1}g)-\mf(f)\right|<\frac{\ep}{2}.
\]
But, as $f(g_j^{-1}g)=|\la\pi(g_j^{-1}g)\xi|\eta\ra|=|\la \pi(g^{-1}g_j)\eta|\xi\ra|$, by (b), there is some $g\in G$ such that $f(g_j^{-1}g)<\ep/2$ for every $j=1,\ldots,n$, hence
\[
\mf(f)=|\mf(f)|\leq\left|\mf(f)-\sum_j t_j f(g_j^{-1}g)\right|+\sum_j t_jf(g_j^{-1}g)<\ep.
\]
This proves that $\pi$ is weakly mixing.\\
$(b)\ \Rightarrow\ (c)$. If $\pi$ contains a nonzero finite-dimensional subrepresentation, then there exists a unit vector $\xi\in\Hh$ whose orbit $\pi(G)\xi$ is relatively compact with respect to the norm topology. Then, let $\ep>0$ be small enough so that $\ep^2/2+\ep<1/2$. There exists a finite set $S\subset G$ such that 
\[
\pi(G)\xi\subset \bigcup_{s\in S}\{\eta\colon \Vert \eta-\pi(s)\xi\Vert<\ep\}.
\]
If $\pi$ satisfies (b), there exists $g\in G$ such that $|\la \pi(g)\xi|\pi(s)\xi\ra|<\ep$ for every $s\in S$. Thus, let $s=s(g)\in S$ be such that 
$\Vert \pi(g)\xi-\pi(s)\xi\Vert<\ep$. Then we would get
\begin{align*}
\ep^2 
&\geq 
\Vert \pi(g)\xi-\pi(s)\xi\Vert^2\\
&= 
2-2\re\la\pi(g)\xi|\pi(s)\xi\ra\\
&\geq 
2(1-\ep),
\end{align*}
which would imply that $1\leq \ep^2/2+\ep<1/2$, a contradiction.\\
$(c)\ \Rightarrow\ (a)$. Assume that $\pi$ is not weakly mixing. Then there exist unit vectors $\xi,\eta\in\Hh$ such that $\mf(|\xi*_\pi\bar\eta|^2)>0$. But
\[
|\xi*_\pi\bar\eta|^2(g)=\la\pi\otimes\bar\pi(g)\xi\otimes\bar\xi|\eta\otimes\bar\eta\ra
\]
where $(\bar\pi,\bar\Hh)$ denotes the conjugate representation of $\pi$ (\cite[Definition A.1.10]{BHV}).
Hence, by formula (\ref{proj}), 
\[
\mf(|\xi*_\pi\bar\eta|^2)=\la P_{\pi\otimes\bar\pi}\xi\otimes\bar\xi|\eta\otimes\bar\eta\ra>0
\]
which shows that the representation $\pi\otimes\bar\pi$ contains a nonzero invariant vector, and by \cite[Proposition A.1.12]{BHV}, that $\pi$ has a nonzero finite-dimensional subrepresentation.
\hfill $\square$

\medskip

The following proposition is a consequence of the fact that every unitary representation is a direct sum of cyclic representations \cite[Proposition C.4.9]{BHV}. 

\begin{prop}\label{Prop2.6}
Let $G$ be a topological group. Then it admits a weakly mixing representation such that $1_G\prec\pi$ if and only if there exists a generalized sequence $(\f_i)$ of normalized functions of positive type on $G$ such that:
\begin{itemize}
\item for every $i$, $\mf(|\f_i|)=0$;
\item for every compact subset $Q\not=\emptyset$ of $G$, one has
\[
\lim_i \Big(\sup_{g\in Q}|\f_i(g)-1|\Big)=0.
\]
\end{itemize}
\end{prop}

For future use, recall that a topological group $G$ is \textit{minimally almost periodic} if its only finite-dimensional unitary representation is the trivial one (see \cite[Section 4.C.d]{BH}). Then we make the following observation. 

\begin{prop}\label{prop2.7}
Let $G$ be a minimally almost periodic group and let $(\pi,\Hh)$ be a unitary representation of $G$. Then $\pi$ is weakly mixing if and only if it has no nonzero invariant vector.
\end{prop}

\section{A weakly mixing approximation property for tracial von Neumann algebras}

Let $(M,\tau)$ be a tracial von Neumann algebra where $\tau$ is a normal, faithful tracial state on $M$.
We denote by $(M)_1$ the unit ball of $M$ with respect to the operator norm, by $U_M$ its unitary group and by $Z_M$ its center. 

We denote by $L^2M$ the usual $M$-bimodule obtained by completion of $M$ for the scalar product $(x,y)\mapsto \tau(y^*x)$, and by $\Vert\cdot\Vert_2$ the corresponding norm, and when it is necessary to specify which trace is used, it is denoted by $\Vert\cdot\Vert_{2,\tau}$. We denote by $x\mapsto \hat x$ the embedding of $M$ into the corresponding dense subspace of $L^2M$. The bimodule structure comes from the natural product $x\hat y z=\wh{xyz}$ for all $x,y,z\in M$.

The group $U_M$ is equipped with the bi-invariant metric $d(u,v)\coloneqq \Vert u-v\Vert_2$ which makes it a complete metric topological group (which is Polish if $M$ has separable predual).

We denote by $\CP(M)$ the convex set of all normal, completely positive maps $\phi:M\rightarrow M$, and we say that $\phi\in \CP(M)$ is \textit{subunital} if moreover $\phi(1)\leq 1$.

The set $\CP_{\tau}(M)$ is the subset of all \textit{subtracial} (i.e. $\tau\circ\phi\leq \tau$) elements of $\CP(M)$. Let $\phi\in \CP_\tau(M)$; the linear map $\hat\phi:\hat x\mapsto \wh{\phi(x)}$ extends to a bounded operator on $L^2M$, and moreover, by \cite[Lemma 1.1.5]{Popa}, there exists a $\phi^*\in \CP(M)$ such that $\tau(\phi(x)y)=\tau(x\phi^*(y))$ for all $x,y\in M$. It is subunital, and if $\phi$ is subunital, then $\phi^*$ belongs to $\CP_\tau(M)$. Finally, one has $\widehat{\phi^*}=\hat\phi^*$.

To every $\phi\in\CP(M)$ corresponds a \textit{pointed} $M$-bimodule $(\Hh_\phi,\xi_\phi)$ characterized by the following properties (see \cite[Subsection 1.1.2]{Popa}) :
\begin{itemize}
\item the vector $\xi_\phi\in\Hh_\phi$ is \textit{cyclic} in the sense that the subspace formed by all finite sums $\sum_i x_i\xi_\phi y_i$, $x_i,y_i\in M$ is dense in $\Hh_\phi$; the latter subspace is denoted by $\Span(M\xi_\phi M)$;
\item the scalar product on $\Hh_\phi$ is entirely determined by the following equality:
\begin{equation}\label{scalprod}
    \la x_1\xi_\phi y_1|x_2\xi_\phi y_2\ra=\tau(\phi(x_2^*x_1)y_1y_2^*)\quad (x_1,x_2,y_1,y_2\in M).
\end{equation}
\end{itemize}
As above, if we need to specify that the particular trace $\tau$ is used with the cp map $\phi$, we will write $\la\cdot|\cdot\ra_\tau$ for the above scalar product.

Let $\Hh$ be a $M$-bimodule. Its \textit{adjoint} $\overline{\Hh}$ is the conjugate Hilbert space of $\Hh$ equiped with the actions of $M$ defined by $x\cdot \ov\xi\cdot y=\ov{y^*\xi x^*}$ for all $x,y\in M$ and $\ov\xi\in \ov{\Hh}$. If $\phi\in \CP_\tau(M)$ then one has $(\ov\Hh,\ov{\xi_\phi})=(\Hh_{\phi^*},\xi_{\phi^*})$.

\medskip

Following \cite[Section 1.1.2]{Popa}, one associates a cp map to the following pairs of pointed $M$-bimodules $(\Hh,\xi)$: let $(\Hh,\xi)$ be a pointed $M$-bimodule such that $\la\xi\cdot|\xi\ra\leq c\tau$ for some $c>0$. Let $T:L^2M\rightarrow \Hh$ be the unique bounded operator defined by $T\hat y=\xi y$ for every $y\in M$. Then $\la\xi y|\xi y\ra\leq c\Vert \hat y\Vert_2^2$, so that $\Vert T\Vert\leq \sqrt{c}$. It turns out that $\phi_{(\Hh,\xi)}\coloneq T^*xT$ belongs to $M$ for every $x\in M$, so that it defines a cp map on $M$. If we denote by $\Hh'$ the closure of $\Span(M\xi M)$, then $U:\Hh_\phi\rightarrow \Hh'$, defined by $U(x\xi_\phi y)=x\xi y$, is easily seen to be an $M$-bimodule isomorphism. 

By functional calculus, it is easy to prove that, given a Hilbert $M$-bimodule $\Hh$, the subset $\Hh^0$ of \textit{left and right bounded vectors} $\xi\in\Hh$ for which there exists a constant $c>0$ such that $\la\xi\cdot |\xi\ra\leq c\tau$ and $\la\cdot\ \xi |\xi\ra\leq c\tau$ is dense in $\Hh$; thus every Hilbert $M$-bimodule is a direct sum of bimodules associated to cp maps of the form $\phi_{(\Hh,\xi)}$ with $\xi\in \Hh^0$. 

\medskip

Let $\Hh$ be a $M$-bimodule. The following unitary representation $(\pi_\Hh,\Hh)$ of $U_M$ plays a central role here; it is defined by
\[
\pi_\Hh(u)\xi=u\xi u^* \quad (u\in U_M,\xi\in \Hh).
\]
For $\phi\in \CP(M)$, we write $\pi_\phi$ instead of $\pi_{\Hh_\phi}$.

\medskip
Motivated by Proposition \ref{Prop2.6}, we propose the following definition.

\begin{defn}\label{defwmap}
Let $(M,\tau)$ be a tracial von Neumann algebra. 
\begin{enumerate}
\item [(1)] Let $\phi\in \CP(M)$. Then we say that it is \textit{weakly mixing} if the representation $\pi_\phi$ is weakly mixing; we denote by $\WM(M)$ the set of all weakly mixing elements of $\CP(M)$, and we also set $\WM_\tau(M)=\WM(M)\cap\CP_\tau(M)$.
\item [(2)] We say that $(M,\tau)$ has the \textit{weakly mixing approximation property} (\textit{WMAP}) if there exists a generalized sequence $(\phi_i)_{i\in I}\subset \WM(M)$ such that
\[
\lim_i \Vert \phi_i(x)-x\Vert_{2,\tau}=0 \quad (x\in M).
\]
\end{enumerate}
\end{defn}

\begin{rem}
A notion of (left) weakly mixing property for $M$-bimodules was introduced in \cite{PS} and was used by H. Tan in \cite{Tan} to prove among others in \cite[Theorem 1.3]{Tan} that, for a II$_1$ factor $M$, the negation of property (T) for $M$ implies the existence of weakly mixing bimodules almost having $M$-central vectors.

Recall that the $M$-bimodule $\Hh$ is \textit{(left) weakly mixing} if it satisfies the following equivalent conditions (see \cite[Theorem A.2.2]{Bou}):
\begin{enumerate}
\item [(i)] there exists a sequence $(u_n)\subset U_M$ such that 
\[
\lim_{n\to\infty}\Big(\sup_{b\in (M)_1}|\la u_n\xi b|\eta\ra|\Big)=0
\]
for all $\xi,\eta\in \Hh_\phi$;
\item [(ii)] the $M$-bimodule $\Hh\otimes_M\ov{\Hh}$ has no nonzero central vector.
\end{enumerate}
Given $\phi\in \CP(M)$, the weakly mixing property of $\phi$ defined here is certainly weaker than the left weakly mixing property of the $M$-bimodule $\Hh_\phi$. Indeed, if $\Hh_\phi$ is (left) weakly mixing, then by condition (i) above, $\pi_\phi$ satisfies condition (b) in Proposition \ref{Prop2.4}, which means that $\phi$ is weakly mixing in our sense. We do not know whether the converse holds true, though.  

Nevertheless, we make the following observation: suppose that $M$ is diffuse; condition (ii) above applied to the $M$-bimodule $\Hh_\phi$ means that the unitary representation $\pi_{\Hh_\phi\otimes_M\ov{\Hh_{\phi}}}$ of $U_M$ by conjugation  has no nonzero invariant vector. But, by \cite[Theorem 4.4]{Jol1}, $U_M$ is minimally almost periodic and thus Proposition 2.7 implies that $\Hh_\phi$ is (left) weakly mixing if and only if $\pi_{\Hh_\phi\otimes_M\ov{\Hh_{\phi}}}$ is a weakly mixing representation.
\end{rem}

Before studying WMAP further, let us prove the following result, which states in particular that, if $M$ is diffuse and if it has the Haagerup property (\cite{Cho,Jol}), then it has the WMAP. Following a suggestion of the referee, we remind the reader that the Haagerup property is characterized in \cite[Theorem 3.4]{BF} and in \cite[Theorem 9]{OOT} by the existence of a \textit{strictly mixing} $M$-bimodule that weakly contains the identity $M$-bimodule $L^2M$; for bimodules, strict mixing is stronger than weak mixing (\cite[Definition 4]{OOT}).

\begin{prop}\label{HAP}
Let $(M,\tau)$ be a tracial von Neumann algebra.
\begin{enumerate}
\item [(1)] If $(M,\tau)$ has the WMAP then so does $(zM,\tau_z)$ for every nonzero central projection $z$.
\item [(2)] If $M$ is diffuse and if $\phi\in \CP_\tau(M)$ is $L^2$-compact, i.e. the operator $\hat\phi\in B(L^2M)$ is compact, then $\phi\in \WM_\tau(M)$. In particular, if $M$ has the Haagerup property then it has the WMAP.
\item [(3)] If $M$ is atomic, then $\WM(M)=\{0\}$, and hence $M$ does not have the WMAP.
\end{enumerate}
\end{prop}
\textsc{Proof.} (1) Recall first that the normalized trace $\tau_z$ on $zM$ is defined by
\[
\tau_z(zx)=\frac{1}{\tau(z)}\tau(zx)\quad (x\in M).
\]
Let then $(\phi_i)_{i\in I}\subset \WM(M)$ be as in Definition \ref{defwmap}. For every $i\in I$, let $\psi_i\in \CP(M)$ be defined by $\psi_i(zx)=z\phi_i(zx)$ for every $x\in M$. Then one has for all $x_1,x_2,y_1,y_2\in M$ and $u\in U_M$:
\begin{align*}
[(zx_1)\xi_{\psi_i}(zy_1)*_{\pi_{\psi_i}}\ov{(zx_2)\xi_{\psi_i}(zy_2)}](zu)
&=
\la zu(zx_1)\xi_{\psi_i}(zy_1)(zu^*)|(zx_2)\xi_{\psi_i}(zy_2)\ra\\
&=
\frac{1}{\tau(z)}\tau(\phi_i(zx_2^*zux_1)zy_1uy_2^*)\\
&=
\frac{1}{\tau(z)}\la u(zx_1)\xi_{\phi_i}(zy_1)u^*|(zx_2)\xi_{\phi_i}(zy_2)\ra\\
&=
\frac{1}{\tau(z)}\cdot[(zx_1)\xi_{\phi_i}(zx_2)*_{\pi_{\phi_i}}\ov{(zx_2)\xi_{\phi_i}(zy_2)}](u).
\end{align*}
This proves that $(\psi_i)_{i\in I}\subset \WM(zM)$. Finally, one has for every $x\in M$:
\[
\Vert \psi_i(zx)-zx\Vert_{2,\tau_z}=\frac{1}{\sqrt{\tau(z)}}\Vert z\phi_i(zx)-zx\Vert_{2,\tau}\leq
\frac{1}{\sqrt{\tau(z)}}\Vert \phi_i(zx)-zx\Vert_{2,\tau}\to_i 0.
\]
(2) Let $x_1,y_1,x_2,y_2\in M$ and $S\subset U_M$ be a finite set. Since $M$ is diffuse, $U_M$ contains a sequence $(u_n)_{n\geq 1}$ which converges weakly to $0$. Using Proposition \ref{Prop2.4}(b), it suffices to prove that 
\begin{equation}\label{convfaible}
    \sum_{s\in S}|\la \pi_\phi(u_ns)x_1\xi_\phi y_1|x_2\xi_\phi y_2\ra|\to 0
\end{equation}
as $n\to\infty$.\\
We have, for $s\in S$ and $n\geq 1$:
\begin{align*}
\la \pi_\phi(u_ns)x_1\xi_\phi y_1|x_2\xi_\phi y_2\ra
&=
\la u_nsx_1\xi_\phi y_1|x_2\xi_\phi y_2u_ns\ra\\
&=
\tau(\phi(x_2^*u_nsx_1)y_1s^*u_n^*y_2^*)\\
&=
\la \hat\phi(\wh{x_2^*u_nsx_1})|\wh{y_2u_nsy_1^*}\ra.
\end{align*}
But $\wh{x_2^*u_nsx_1}\to 0$ weakly as $n\to \infty$, thus, as $\hat\phi$ is a compact operator, one gets
\[
\Vert \hat\phi(\wh{x_2^*u_nsx_1})\Vert_2\to_{n\to\infty}0.
\]
Finally,
\[
\sum_{s\in S}|\la \pi_\phi(u_ns)x_1\xi_\phi y_1|x_2\xi_\phi y_2\ra|
\leq \sum_{s\in S}\Vert \hat\phi(\wh{x_2^*u_nsx_1})\Vert_2\Vert y_2\Vert\Vert y_1\Vert_2
\to 0
\]
as $n\to\infty$, which ends the proof of (\ref{convfaible}).\\
(3) Assume that $M$ is a finite, atomic von Neumann algebra. As is well-known, it admits a direct sum decomposition 
\[
M=\bigoplus_{i\in I} N_i
\]
where $N_i$ is a matrix algebra $M_{n_i}(\C)$ with $n_i\geq 1$ for every $i$. If $I$ was uncountable, then $M$ would admit no faithful, normal tracial state, thus $I$ is at most countable. The faithful, normal tracial state $\tau$ on $M$ is of the form
\[
\tau(\oplus_i x_i)=\sum_i \mu_i \Tr_{n_i}(x_i)\quad (\oplus_i x_i\in M)
\]
where $\mu_i>0$ for every $i\in I$ and $\sum_i \mu_i n_i=1$. The topology on $U_M$ is metrisable with respect to $\Vert\cdot\Vert_2$ and $U_M=\prod_i U(n_i)$ is a compact group. This implies that no nonzero element $\phi\in \CP(M)$ is weakly mixing, as all non trivial unitary representations of $U_M$ are direct sums of finite-dimensional ones.
\hfill $\square$

\medskip
One of our next tasks is to prove that the generalized sequence $(\phi_i)_{i\in I}$ in Definition \ref{defwmap}(2) can be chosen in $\WM_\tau(M)$. We will apply techniques borrowed from the proof of Theorem 2.2 of \cite{BF}.

First, the next lemma displays stability properties of suitably modified weakly mixing cp maps.

\begin{lem}\label{lem3.2}
Let $(M,\tau)$ be a tracial von Neumann algebra.
\begin{enumerate}
\item [(1)] $\WM(M)$ is a subcone of $\CP(M)$: for all $\lm\geq 0, \phi,\psi\in \WM(M)$, one has $\lm\phi+\psi\in \WM(M)$.
\item [(2)] Let $\phi\in \WM(M)$ and $a\in (M)_1$. Define $\psi,\rho:M\rightarrow M$ by $\psi(x)=a^*\phi(x)a$ and $\rho(x)=\phi(a^*xa)$ for every $x\in M$. Then $\psi$ and $\rho$ belong to $\WM(M)$ and more precisely, one has the following equalities for all $x_1,x_2,y_1,y_2\in M$:
\begin{equation}\label{coeff1}
(x_1\xi_\psi y_1)*_{\pi_\psi}(\ov{x_2\xi_\psi y_2})=(x_1\xi_\phi ay_1)*_{\pi_\phi}(\ov{x_2\xi_\phi ay_2})   
\end{equation} 
and 
\begin{equation}\label{coeff2}
(x_1\xi_\rho y_1)*_{\pi_\rho}(\ov{x_2\xi_\rho y_2})=(x_1a\xi_\phi y_1)*_{\pi_\phi}(\ov{x_2a\xi_\phi y_2}).
\end{equation}
Moreover, if $\phi\in\WM_\tau(M)$ then $\psi,\rho\in \WM_\tau(M)$.
\item [(3)] Let $z\in Z_M$ be a nonzero projection and let $\theta\in \WM_{\tau_z}(zM)$, where $\tau_z$ is as in Proposition \ref{HAP}. Set $\phi(x)=\theta(zx)=z\theta(zx)$ for $x\in M$. Then $\phi\in \WM_\tau(M)$, and we have for all $x_1,x_2,y_1,y_2\in M$ and $u\in U_M$:
\begin{equation}\label{coeff3}
(x_1\xi_\phi y_1)*_{\pi_\phi}(\ov{x_2\xi_\phi y_2})(u)=\tau(z)(zx_1\xi_\theta zy_1)*_{\pi_\theta}(\ov{zx_2\xi_\theta zy_2})(zu).
\end{equation}
\end{enumerate}
\end{lem}
\textsc{Proof.} (1) It is obvious that $\lm\phi\in \WM(M)$ if $\phi\in\WM(M)$. Next, it is easily seen that, for all $x_1,x_2,y_1,y_2\in M$ one has
\[
(x_1\xi_{\phi+\psi}y_1)*_{\pi_{\phi+\psi}}(\ov{x_2\xi_{\phi+\psi}y_2})=(x_1\xi_\phi y_1)*_{\pi_\phi}(\ov{x_2\xi_\phi y_2})
+(x_1\xi_\psi y_1)*_{\pi_\psi}(\ov{x_2\xi_\phi y_2}),
\]
which implies that 
\begin{multline*}
\mf(|(x_1\xi_{\phi+\psi}y_1)*_{\pi_{\phi+\psi}}(\ov{x_2\xi_{\phi+\psi}y_2})|)\\
\leq \mf(|(x_1\xi_\phi y_1)*_{\pi_\phi}(\ov{x_2\xi_\phi y_2})|)+\mf(|(x_1\xi_\psi y_1)*_{\pi_\psi}(\ov{x_2\xi_\phi y_2})|)=0.
\end{multline*}
(2) The maps $\psi$  and $\rho$ are obviously cp.
Let us establish equality (\ref{coeff1}), which will prove that $\psi\in\WM(M)$: one has for every $u\in U_M$,
\begin{align*}
(x_1\xi_\psi y_1)*_{\pi_\psi}(\ov{x_2\xi_\psi y_2})(u)
&=
\la ux_1\xi_\psi y_1u^*|x_2\xi_\psi y_2\ra\\
&=
\tau(\psi(x_2^*ux_1)y_1u^*y_2^*)=\tau(\phi(x_2^*ux_1)ay_1u^*(ay_2)^*)\\
&=
\la ux_1\xi_\phi ay_1u^*|x_2\xi_\phi ay_2\ra=
(x_1\xi_\phi ay_1)*_{\pi_\phi}(\ov{x_2\xi_\phi ay_2})(u).
\end{align*}
Equality (\ref{coeff2}) is verified similarly.

If moreover $\phi\in \WM_\tau(M)$ then one has for every $x\in M$
\[
\tau(\psi(x^*x))=\tau(a^*\phi(x^*x)a)=\tau(\phi(x^*x)^{1/2}aa^*\phi(x^*x)^{1/2})\leq\tau(\phi(x^*x)\leq \tau(x^*x),
\]
which proves that $\psi\in \CP_\tau(M)$, and similarly for $\rho$.\\
(3) The verification of the fact that $\phi\in \WM_\tau(M)$ and equality (\ref{coeff3}) follow from part (2). 
\hfill $\square$

\medskip
The following result is strongly influenced by \cite[Theorem 22]{BF}. We present a detailed proof for the sake of completeness. We use Lemmas 2.3 and 2.4 of \cite{BF}, which are parts of \cite[Lemma 1.1.2]{Popa}.

\begin{lem}\label{lemBF}
(\cite[Lemmas 2.3 and 2.4]{BF})
Let $(M,\tau)$ be a tracial von Neumanna algebra.
\begin{enumerate}
\item [(1)]  Let $\f\in \WM(M)$. Set $a=1\vee\f(1)$, and define $\psi\in \CP(M)$ by $\psi(x)=a^{-1/2}\f(x)a^{-1/2}$ for every $x\in M$. Then $\psi\in \WM(M)$, it is subunital, and it has the following properties: $\tau\circ\psi\leq \tau\circ\f$ and
\[
\Vert \psi(x)-x\Vert_2\leq \Vert\f(x)-x\Vert_2+2\Vert \f(1)-1\Vert_2^{1/2}\quad (x\in (M)_1).
\]
\item [(2)] Let $\psi\in \WM(M)$ be subunital, and set 
\[
b=1\vee\Big(\frac{d\tau\circ\psi}{d\tau}\Big),
\] 
and define $\phi\in \CP(M)$ by $\phi(x)=\psi(b^{-1/2}xb^{-1/2})$ for every $x\in M$. Then $\phi\in \WM_\tau(M)$, and it satisfies the following inequalities:
$\phi(1)\leq \psi(1)\leq 1$ and 
\[
\Vert \phi(x)-x\Vert_2^2\leq 2\Vert \psi(x)-x\Vert_2+5\Vert \tau\circ\psi-\tau\Vert^{1/2}\quad (x\in (M)_1).
\]
\end{enumerate}
\end{lem}
\textsc{Proof.} (1) As $a\geq 1$, we have $0\leq a^{-1}\leq 1$, and, for every $x\in (M)_1$, we get
\begin{align*}
\Vert \psi(x)-x\Vert_2 &\leq
\Vert a^{-1/2}(\f(x)-x)a^{-1/2}\Vert_2+\Vert a^{-1/2}x(a^{-1/2}-1)\Vert_2+\Vert (a^{-1/2}-1)x\Vert_2\\
& \leq
\Vert \f(x)-x\Vert_2+2\Vert a^{-1/2}-1\Vert_2\leq \Vert \f(x)-x\Vert_2+2\Vert a^{-1}-1\Vert_1^{1/2}\\
& \leq
\Vert \f(x)-x\Vert_2+2\Vert a^{-1}\Vert \Vert 1-a\Vert_1^{1/2}\leq 
\Vert \f(x)-x\Vert_2+2\Vert 1-a\Vert_2^{1/2}\\
&\leq 
\Vert \f(x)-x\Vert_2+2\Vert \f(1)-1\Vert_2^{1/2},
\end{align*}
where we used Powers-St\o rmer inequality in the second line, and $\Vert 1-a\Vert_1\leq \Vert 1-a\Vert_2$ by Cauchy-Schwarz inequality.
\par\noindent
(2) We divide this part into two steps.\\
\textbf{Step 1} The following inequalities hold true:
\[
\tau(b)\leq 2,\quad \Vert y\Vert_2^2\leq \Vert y\Vert\Vert y\Vert_1\quad \textrm{and} \quad \Vert \phi(y)\Vert_1\leq \Vert y\Vert_1\quad (y\in M).
\]
Indeed, one has $\tau(b)\leq \tau(1+\psi(1))\leq 2$ since $\psi(1)\leq 1$, and $\Vert y\Vert_2^2=\tau(|y|^2)\leq \Vert y\Vert \tau(|y|)=\Vert y\Vert\Vert y\Vert_1$. In order to prove the last inequality, recall that, as $\phi$ is subtracial, there exists a unique $\phi^*\in \CP(M)$ such that $\tau(x\phi(y))=\tau(\phi^*(x)y)$ for all $x,y\in M$. It is straightforward to check that furthermore $\phi^*(1)\leq 1$. Then we have for every $y\in M$:
\begin{multline*}
\Vert \phi(y)\Vert_1
=
\sup\{|\tau(\phi(y)z)|\colon z\in (M)_1\}\\
=\sup\{|\tau(y\phi^*(z))|\colon z\in (M)_1\}
\leq
\sup\{|\tau(yz)|\colon z\in (M)_1\}=\Vert y\Vert_1.
\end{multline*}
\textbf{Step 2} Fix $x\in (M)_1$. Then, by Step 1
\begin{align*}
\Vert \phi(x)-x\Vert_2^2
& \leq 
\Vert \phi(x)-x\Vert\Vert \phi(x)-x\Vert_1\leq 2\Vert \phi(x)-x\Vert_1
 \leq
2\Vert \phi(x)-\underbrace{\phi(b^{1/2}xb^{1/2})}_{=\psi(x)}\Vert_1\\
& +2\underbrace{\Vert \psi(x)-x\Vert_1}_{\leq \Vert \psi(x)-x\Vert_2}
 \leq
2\Vert x-b^{1/2}xb^{1/2}\Vert_1+2\Vert \psi(x)-x\Vert_2\\
& \leq
2(\Vert x(1-b^{1/2})\Vert_1+\Vert (1-b^{1/2})xb^{1/2}\Vert_1+\Vert \psi(x)-x\Vert_2)\\
& \leq
2(\Vert x\Vert_2\Vert 1-b^{1/2}\Vert_2+\Vert 1-b^{1/2}\Vert_2\underbrace{\Vert xb^{1/2}\Vert_2}_{\leq \tau(b)^{1/2}\leq \sqrt 2}+\Vert \psi(x)-x\Vert_2)\\
& \leq
2\Vert \psi(x)-x\Vert_2+(2+2\sqrt 2)\Vert 1-b^{1/2}\Vert_2
 \leq 
2\Vert \psi(x)-x\Vert_2+5\Vert 1-b\Vert_1^{1/2}\\
& \leq 
2\Vert \psi(x)-x\Vert_2+
5\Vert \tau\circ\psi-\tau\Vert^{1/2},
\end{align*}
where we used Powers-St\o rmer inequality once again.
\hfill $\square$

\begin{thm}\label{subunital}
Let $(M,\tau)$ be a tracial von Neumann algebra. The following conditions are equivalent.
\begin{enumerate}
\item [(1)] $(M,\tau)$ has the WMAP.
\item [(2)] There exists a generalized sequence $(\phi_j)_{j\in \J}\subset \WM_\tau(M)$ such that
$\phi_j(1)\leq 1$ for every $j$ and
\[
\lim_j\Vert \phi_j(x)-x\Vert_2=0
\]
for every $x\in M$.
\item [(3)] There exists a generalized sequence of pointed $M$-bimodules $((\Hh_j,\xi_j))_{j\in \J}$ with the following properties:
\begin{enumerate}
\item $\displaystyle{\lim_j(\Vert \la\cdot\ \xi_j|\xi_j\ra-\tau\Vert+\Vert \la\xi_j\cdot |\xi_j\ra-\tau\Vert)=0}$;
\item $\displaystyle{\lim_j \Vert x\xi_j-\xi_j x\Vert=0}$ for every $x\in M$;
\item the unitary representation $\pi_{\Hh_j}$ is weakly mixing for every $j\in\J$.
\end{enumerate}
\end{enumerate}
\end{thm}
\textsc{Proof.} $(1)\ \Rightarrow\ (2)$.
Let $(\f_i)_{i\in I}\subset \WM(M)$ be a generalized sequence as in part (2) of Definition \ref{defwmap}. Set $a_i=1\vee \f_i(1)$ and $\psi_i(x)=a_i^{-1/2}\f_i(x)a_i^{-1/2}$ for all $x\in M$ and $i\in I$. Then, by part (1) of Lemma \ref{lemBF}, each $\psi_i$ is subunital and we have 
\begin{equation}\label{limpsii}
\lim_i\Vert \psi_i(x)-x\Vert_2=0\quad \textrm{and}\quad \lim_i \tau\circ\psi_i(x)=\tau(x) \quad (x\in M).
\end{equation} 
Let $K$ be the weak closure of the convex hull of $(\tau\circ\psi_i)_{i\in I}$ in $M_*$. As $K$ is the norm closure of the above mentioned convex hull, there exists a net $(\psi'_j)_{j\in \J}$ of convex combinations of the $\psi_i$'s such that 
\[
\lim_j\Vert \tau\circ\psi'_j-\tau\Vert=0\quad \textrm{and}\quad \lim_j\Vert\psi'_j(x)-x\Vert_2=0 \quad (x\in M).
\] 
Moreover, the $\psi'_j$'s are subunital and belong to $\WM(M)$ by Lemma \ref{lem3.2}(1).

For every $j\in \J$, set 
\[
b_j=1\vee\Big(\frac{d\tau\circ\psi'_j}{d\tau}\Big)
\]
and define $\phi_j\in \CP(M)$ by $\phi_j(x)=\psi'_j(b_j^{-1/2}xb_j^{-1/2})$ for $x\in M$. Then, by part (2) of Lemma \ref{lemBF}, every $\phi_j$ is subunital, $(\phi_j)_{j\in\J}\subset \WM_\tau(M)$ and $\Vert \phi_j(x)-x\Vert_2\to 0$ for every $x\in M$.\\
$(2)\ \Rightarrow\ (3)$. Let $(\phi_j)_{j\in \J}\subset WM_\tau(M)$ be as in (2). Set $(\Hh_j,\xi_j)=(\Hh_{\phi_j},\xi_{\phi_j})$ for every $j\in \J$.

Let $x\in (M)_1$. One has
\[
|\la x\xi_j|\xi_j\ra-\tau(x)|=|\tau(\phi_j(x)-x)|=|\tau(x(\phi_j^*(1)-1))|\leq \Vert x\Vert_2\Vert \phi_j^*(1)-1\Vert_2\leq 
\Vert \phi_j^*(1)-1\Vert_2
\]
and, as $\phi_j^*$ is subunital,
\[
\Vert\phi_j^*(1)-1\Vert_2^2
=
\tau(\phi_j^*(1)^2-2\phi_j^*(1)+1)\leq \tau(1-\phi_j^*(1))=\tau(1-\phi_j(1))
\leq \Vert 1-\phi_j(1)\Vert_2\to_j 0.
\]
Similarly, 
\[
|\la \xi_j x|\xi_j\ra-\tau(x)|=|\tau(\phi_j(1)x-x)|\leq \Vert x\Vert_2\Vert\phi_j(1)-1\Vert_2\leq \Vert\phi_j(1)-1\Vert_2\to_j0.
\]
This proves (a), and (c) is obvious. In order to prove (b), we observe first that, for all $j\in\J$ and $x\in M$,
\[
\Vert \xi_j x\Vert_2^2=\la \xi_{\phi_i}x|\xi_{\phi_i}x\ra=\tau(\phi_j(1)xx^*)\leq \Vert x\Vert_2^2
\]
and similarly $\Vert x\xi_j\Vert^2=\tau(\phi_j(x^*x))\leq \Vert x\Vert_2^2$. Hence, one gets for every $x\in M$
\begin{multline*}
\Vert x\xi_j-\xi_j x\Vert^2=\Vert x\xi_j\Vert^2+\Vert \xi_j x\Vert^2-2\re\la x\xi_j|\xi_j x\ra\\
\leq 2\Vert x\Vert_2^2-2\re\tau(\phi_j(x)x^*)=2\re\tau((x-\phi_j(x)x^*)\leq 2\Vert x-\phi_j(x)\Vert_2\Vert x\Vert_2\to_j0
\end{multline*}
which proves (b).\\
$(3)\ \Rightarrow\ (1)$. By density of the set $\Hh_j^0$ of left and right bounded vectors in $\Hh_j$, we assume that every element of the generalized sequence $(\xi_j)_{j\in \J}$ is left and right bounded, and we denote by $T_j:L^2M\rightarrow \Hh_j$ the operator $T_j\hat y=\xi_j y$ for every $y\in M$. We set then $\phi_j(x)=T_j^*xT_j$ for every $x\in M$, so that $\phi_j\in \WM(M)$. One has for every $x\in M$ and every $j\in\J$:
\[
\la x\xi_j|\xi_j x\ra=\la x\xi_j|T_j\hat{x}\ra=\la \xi_j|x^*T_j\hat{x}\ra=\la \hat 1|T_j^*x^*T_j\hat{x}\ra=\tau(\phi_j(x)x^*).
\]
Hence, for fixed $x\in (M)_1$,
\begin{align*}
\Vert \phi_j(x)-x\Vert_2^2
&=
2\Vert x\Vert_2^2-2\re\tau(\phi_j(x)x^*)=2\Vert x\Vert_2^2-\re\la x\xi_j|\xi_j x\ra\\
&=
(\tau(x^*x)- \la x^*x\xi_j|\xi_j\ra)+(\tau(xx^*)-\la \xi_j xx^*|\xi_j\ra)\\
&+
\Vert x\xi_j\Vert^2+\Vert \xi_j x\Vert^2-2\re\la x\xi_j|\xi_j x\ra\\
&\leq 
\Vert \la \cdot\ \xi_j|\xi_j\ra-\tau\Vert+\Vert \la \xi_j\cdot|\xi_j\ra-\tau\Vert+\Vert x\xi_j-\xi_jx\Vert^2\to_j0
\end{align*}
This shows that $(M,\tau)$ has the WMAP.
\hfill $\square$

\medskip
We prove now that the WMAP does not depend on the chosen trace.

\begin{prop}\label{central}
Let $(M,\tau)$ be a tracial von Neumann algebra which has the WMAP, and let $\tau'$ be a normal, normalized, faithful trace on $M$. Then $(M,\tau')$ has the WMAP.
\end{prop}
\textsc{Proof.}  Let us prove first the following statement: 

\textit{Let $(z_n)_{n\geq 1}\subset Z_M$ be an increasing sequence of nonzero projections which converges strongly to $1$. Let $\tau_n$ be the normalized trace $z_nx\mapsto \tau(z_nx)/\tau(z_n)$ on $z_nM$. If each pair $(z_nM,\tau_n)$ has the WMAP, then so does $(M,\tau)$.}\\
Indeed, let $F\subset M$ be a finite set and let $\ep>0$. Choose first $n$ large enough so that $\Vert (z_n-1)x\Vert_{2,\tau}<\ep/2$ for every $x\in F$. Next, choose $\theta\in \WM_{\tau_n}(z_nM)$ so that 
\[
\max_{x\in F}\Vert \theta(z_nx)-z_nx\Vert_{2,\tau_n}<\ep/2.
\]
Define $\phi\in \WM_\tau(M)$ as in Lemma \ref{lem3.2}(2), and observe that, by definition, $\phi(x)=\phi(z_nx)$ for every $x$. We get for every $x\in F$
\begin{align*}
\Vert \phi(x)-x\Vert_{2,\tau}
&\leq
\Vert \phi(z_nx)-z_nx\Vert_{2,\tau}+\Vert (z_n-1)x\Vert_{2,\tau}\\
&=
\sqrt{\tau(z_n)}\Vert \theta(z_nx)-z_nx\Vert_{2,\tau_n}+\Vert(z_n-1)x\Vert_{2,\tau}<\ep,
\end{align*}
where we have used the equality $\Vert zy\Vert_{2,\tau}=\sqrt{\tau(z)}\Vert zy\Vert_{2,\tau_z}$ which holds for any nonzero projection $z\in Z_M$ and every $y\in M$.
This proves that $(M,\tau)$ has the WMAP.

Let now $\tau'$ be another trace on $M$.
There exists a positive, invertible operator $h$ affiliated to $Z_M$ such that $\tau'(x)=\tau(hx)$ for every $x\in M$. 

Let us assume first that $h$ and $h^{-1}$ are bounded, and let $(\phi_i)_{i\in I}\subset \WM_\tau(M)$ be such that 
\[
\lim_i \Vert \phi_i(x)-x\Vert_{2,\tau}=0 \quad (x\in M).
\]
Set for all $i\in I$ and $x\in M$: $\psi_i(x)\coloneqq h^{-1}\phi_i(hx)$. Then $\psi_i$ is cp and, for every $x\in M$,
\[
\tau'\circ\psi_i(x^*x)=\tau(\phi_i(hx^*x))\leq \tau(hx^*x)=\tau'(x^*x),
\]
which shows that $\psi_i\in \CP_{\tau'}(M)$ as well. Let us prove next that $\psi_i\in \WM_{\tau'}(M)$; as in Lemma \ref{lem3.2}, it is easy to check that one has for all $x_1,x_2,y_1,y_2\in M$, where the scalar product on the left hand side is $\la\cdot|\cdot\ra_{\tau'}$:
\[
(x_1\xi_{\psi_i}y_1)*_{\pi_{\psi_i}}(\ov{x_2\xi_{\psi_i}y_2})=(x_1\xi_{\phi_i}y_1)*_{\pi_{\phi_i}}(\ov{hx_2\xi_{\phi_i}y_2}),
\]
which implies that $\psi_i\in\WM_{\tau'}(M)$.\\
Next, for $x\in M$, 
\begin{align*}
\Vert \psi_i(x)-x\Vert^2_{2,\tau'}
&=
\tau'((\psi_i(x^*)-x^*)(\psi_i(x)-x))\\
&=
\tau(h(h^{-1}(\phi_i(hx^*)-x^*)(h^{-1}\phi_i(hx)-x))\\
&=
\tau((\phi_i(hx^*)-hx^*)h^{-1}(\phi_i(hx)-hx))=\tau(h^{-1}|\phi_i(hx)-hx|^2)\\
&\leq
\Vert h^{-1}\Vert\Vert \phi_i(hx)-hx\Vert^2_{2,\tau}\to 0
\end{align*}
as $i\to\infty$.

If $h$ or $h^{-1}$ is not bounded, let $z_n$ be the spectral projection of $h$ corresponding to the interval $[1/n,n]$. Then $h_n=z_nh$ and its inverse are bounded operators of the center of $z_nM$, and the pair $(z_nM,\tau'_n)$ has the WMAP for every $n$. By the first part of the proof, $(M,\tau')$ has the WMAP.
\hfill $\square$

\medskip
As for the Haagerup property for discrete groups, one has:

\begin{thm}
Let $G$ be a discrete group. Then $G$ does not have Kazhdan's property (T) if and only if its von Neumann algebra $L(G)$ has the WMAP.
\end{thm}
\textsc{Proof.} Assume that $G$ does not have Kazhdan's property (T), and let $(\f_i)$ be a generalized sequence of normalized functions of positive type as in Proposition \ref{Prop2.6}. For each $i$, let $m_{\f_i}$ be the multiplier on $L(G)$ characterized by
\[
m_{\f_i}(\lm(g))=\f_i(g)\lm(g)\quad (g\in G),
\]
where $\lm$ denotes the left regular representation whose range generates $L(G)$ (see for instance \cite[Lemma 2]{Cho}). Then $m_{\f_i}\in \CP_\tau(L(G))$ for every $i$, and even $\tau\circ m_{\f_i}=\tau$. 
By standard approximation, one has $\Vert m_{\f_i}(x)-x\Vert_2\to 0$ for every $x\in L(G)$; see \cite[Theorem 3]{Cho}. 

Thus, we just need to prove that $m_{\f_i}\in \WM_\tau(L(G))$ for every fixed $i$. Set $\phi=m_{\f_i}$, $\xi=\xi_\phi$ (which is a unit vector) and $\pi=\pi_\phi$ for short.
By Proposition \ref{Prop2.4}(b'), all we need to prove is that, for all $s,t,r,w\in G$, for every finite set $F\subset U_{L(G)}$ and for every $\ep>0$, there exists $g\in G$ such that 
\[
\max_{u\in F}|\la\pi(\lm(g)u)\lm(s)\xi\lm(t)|\lm(r)\xi \lm(w)\ra|<\ep.
\]
For every $u\in F$, let us choose $v_u\in L(G)$ which has finite support and such that $\Vert v_u\Vert_2\leq 1$ and $\Vert u-v_u\Vert_2<\ep/4$; let then $K=K^{-1}\subset G$ be a finite, symmetric subset such that $v_u(h)=0$ for all $h\notin K$ and $u\in F$. 

For $g\in G$ and $u\in F$ we have
\begin{align*}
|\la\pi(\lm(g)u)\lm(s)\xi\lm(t)|\lm(r)\xi \lm(w)\ra|
\leq {} &
|\la\lm(g)(u-v_u)\lm(s)\xi\lm(t)|\lm(r)\xi\lm(wg)u\ra|\\
& +
|\la\lm(g)v_u\lm(s)\xi\lm(t)|\lm(r)\xi\lm(wg)v_u\ra|\\
& +
|\la\lm(g)v_u\lm(s)\xi\lm(t)|\lm(r)\xi\lm(wg)(u-v_u)\ra|\\
<{} & 
\ep/2+|\la\lm(g)v_u\lm(s)\xi\lm(t)|\lm(r)\xi\lm(wg)v_u\ra|.
\end{align*}
But
\begin{align*}
\la\lm(g)v_u\lm(s)\xi\lm(t)|\lm(r)\xi\lm(wg)v_u\ra
&=
\tau(\phi[\lm(r^{-1}g)v_u\lm(s)]\lm(tg^{-1}w^{-1})v_u^*)\\
&=
\sum_{h,k\in K}v_u(h)\ov{v_u(k)}\tau(\phi[\lm(r^{-1}ghs)]\lm(tg^{-1}w^{-1}k^{-1}))\\
&=
\sum_{h,k\in K}v_u(h)\ov{v_u(k)}\f_i(r^{-1}ghs)\tau(\lm(r^{-1}ghstg^{-1}w^{-1}k^{-1}))\\
&=
\sum_{h\in K} v_u(h)\ov{v_u(r^{-1}ghstg^{-1}w^{-1})}\f_i(r^{-1}ghs)
\end{align*}
for every $g\in G$. As $\mf(|\f_i|)=0$, by Proposition \ref{Prop2.4}, there exists $g\in G$ such that 
\[
\max_{h\in K}|\f_i(r^{-1}ghs)|<\frac{\ep}{2|K|},
\]
which shows that 
\[
\max_{u\in F}|\la\pi(\lm(g)u)\lm(s)\xi\lm(t)|\lm(r)\xi \lm(w)|<\ep.
\]
Conversely, suppose that $M\coloneqq L(G)$ has the WMAP and let $(\phi_i)\subset \WM(M)$ be as in Definition \ref{defwmap}. For every $i$, define $\f_i:G\rightarrow \C$ by 
\[
\f_i(g)=\tau(\phi_i(\lm(g))\lm(g^{-1}))\quad (g\in G).
\]
Then, by \cite[Lemma 1]{Cho}, $\f_i$ is a function of positive type on $G$, and it is straightforward to check that $\f_i(g)=\la\pi_{\phi_i}(\lm(g))\xi_{\phi_i}|\xi_{\phi_i}\ra$ for every $g\in G$, which proves that $\mf(|\f_i|)=0$ for every $i$. Finally, one obviously has $\lim_i \f_i(g)=1$ for every $g\in G$, which proves that $G$ does not have property (T).
\hfill $\square$

\medskip
Let us recall property (T) for tracial von Neumann algebras from \cite{Popa} and \cite{PP}.

\begin{defn}\label{propT}
A tracial von Neumann algebra $(M,\tau)$ has \textit{property (T)} if it satisfies one (hence both) of the following equivalent conditions:
\begin{enumerate}
\item [(1)] for every $\ep>0$, there exists a finite set $F=F(\ep)\subset M$ and $0<\delta=\delta(\ep)\leq \ep$ such that, if $\phi\in \CP_\tau(M)$ is subunital and if
\[
\max_{x\in F}\Vert \phi(x)-x\Vert_2\leq \delta,
\]
then
\[
\sup_{x\in (M)_1}\Vert \phi(x)-x\Vert_2\leq\ep.
\]
\item [(2)] there exists $F_0\subset M$ finite and $\delta_0>0$ such that for every $\ep>0$, there exists $0<\delta(\ep)\leq \ep$ so that if $\Hh$ is a Hilbert $M$-bimodule with a vector $\xi\in\Hh$ satisfying $\Vert y\xi-\xi y\Vert\leq \delta_0$ for every $y\in F_0$, $\Vert \la\cdot\ \xi|\xi\ra-\tau\Vert<\delta(\ep)$ and  $\Vert \la\xi\cdot |\xi\ra-\tau\Vert<\delta(\ep)$, then there exists $\xi_0\in\Hh$ such that $x\xi_0=\xi_0x$ for every $x\in \Hh$ and $\Vert \la\cdot\ \xi_0|\xi_0\ra-\tau\Vert<\ep$ and $\Vert \la\xi_0\cdot |\xi_0\ra-\tau\Vert<\ep$. (Such a vector $\xi_0$ is called a \textit{central vector}).
\end{enumerate}
 
\end{defn}

It turns out that property (T) for $(M,\tau)$ does not depend on the trace $\tau$: see for instance \cite[Theorem 1]{PP}. Moreover, if $M$ is a II$_1$ factor, the control on the nonzero central vector $\xi_0$ is unnecessary. Indeed, due to \cite[Proposition 1]{CoJo} which proves that condition (wT) below implies (sT), the following three conditions on the II$_1$ factor $M$ are equivalent:
\begin{enumerate}
\item [(wT)] $M$ has property (T): there exists $F\subset M$ finite and $\ep>0$ such that, for every Hilbert $M$-bimodule $\Hh$ and unit vector $\xi\in \Hh$ with
\[
\max_{x\in F}\Vert x\xi-\xi x\Vert<\ep
\]
there is a unit vector $\eta\in \Hh$ which is central.
\item [(mT)] There exists $F_0\subset M$ finite such that, for every $0<\ep<1$, there exists $0<\delta(\ep)\leq \ep$ such that, if $\Hh$ is a Hilbert $M$-bimodule with a unit vector $\xi\in \Hh$ satisfying 
\[
\max_{y\in F_0}\Vert y\xi-\xi y\Vert<\delta(\ep)
\]
and 
\[
\Vert \la\cdot\ \xi|\xi\ra-\tau\Vert+\Vert \la \xi\cdot|\xi\ra-\tau\Vert <\ep
\]
then there exists a unit vector $\xi_0\in \Hh$ which is central.
\item [(sT)] There exists $F_1\subset M$ finite, $\ep_1>0$ and $K>0$ such that, for every $0<\delta\leq \ep_1$, for every pointed $M$-bimodule $(\Hh,\xi)$ where $\xi$ is a unit vector such that
\[
\max_{x\in F_1}\Vert x\xi-\xi x\Vert<\delta
\]
then there is a unit vector $\eta\in \Hh$, which is central and such that
$\Vert \eta-\xi\Vert< K\delta$.
\end{enumerate}

\medskip
Here is our last main result.

\begin{thm}\label{final}
Let $(M,\tau)$ be a tracial von Neumann algebra.
\begin{enumerate}
\item [(1)] If $M$ has the WMAP, then it does not have property (T).
\item [(2)] Assume that $M$ is a $\mathrm{II}_1$ factor. If it does not have property (T) then it has the WMAP.
\end{enumerate}
\end{thm}
\textsc{Proof.} (1) Assume that $M$ has property (T) and take $\ep=1/2$ in Definition \ref{propT}(1); let $F\subset M$ and $\delta>0$ be the associated data. If $\phi\in \CP_\tau(M)$ is subunital and such that
\[
\max_{x\in F}\Vert \phi(x)-x\Vert_2\leq \delta,
\]
then we have in particular
\[
\sup_{u\in U_M}\Vert \phi(u)-u\Vert_2\leq 1/2.
\]
We get for every $u\in U_M$
\[
|\xi_\phi*_{\pi_\phi}\ov{\xi_\phi}(u)-1|=|\tau(\phi(u)u^*)-1)|=|\tau((\phi(u)-u)u^*)|\leq \Vert \phi(u)-u\Vert_2\leq 1/2.
\]
This shows that $\Vert \xi_\phi*_{\pi_\phi}\ov{\xi_\phi}-1\Vert_\infty\leq 1/2$. Hence,
\[
\Vert |\xi_\phi*_{\pi_\phi}\ov{\xi_\phi}|-1\Vert_\infty\leq\Vert \xi_\phi*_{\pi_\phi}\ov{\xi_\phi}-1\Vert_\infty\leq 1/2
\]
and this implies that $|\mf(|\xi_\phi*_{\pi_\phi}\ov{\xi_\phi}|-1|\leq \Vert \xi_\phi*_{\pi_\phi}\ov{\xi_\phi}-1\Vert_\infty\leq 1/2$ and thus $\phi$ cannot belong to $WM_\tau(M)$.\\
(2) Let $M$ be a II$_1$ factor that does not have property (T), hence $M$ satisfies the negation of condition (mT) above. Consider then the set $\J$ of pairs $(F,n)$ where $1\in F\subset M$ is finite and $n\geq 1$ is an integer. We define the following partial order on $\J$: $(F,n)\leq (F',n')$ if and only if $F\subset F'$ and $n\leq n'$. Thus $\J$ is a directed set. As $M$ does not have property (T), for every $j=(F,n)\in\J$, there exists a pointed Hilbert $M$-bimodule $(\Hh_j,\xi_j)$ such that 
\[
\max_{x\in F}\Vert x\xi_j-\xi_j x\Vert \leq 1/n
\]
and 
\[
\Vert \la\cdot\ \xi_j|\xi_j\ra-\tau\Vert+\Vert\la \xi_j\cdot |\xi_j\ra-\tau\Vert\leq 1/n
\]
but such that $\Hh_j$ has no nonzero central vector.  
As $M$ is diffuse, it follows from \cite[Theorem 4.4]{Jol1} that its unitary group is minimally almost periodic, which means by Proposition 2.7 that $\pi_{\Hh_j}$ is weakly mixing. By Theorem \ref{subunital}, $M$ has the WMAP. 
\hfill $\square$

\begin{rem}
Let $(M,\tau)$ be a tracial, atomic von Neumann algebra. It follows from Proposition \ref{HAP} that $M$ does not have the WMAP. We claim that it has property (T).
Indeed, recall from the proof of the above mentioned proposition that $U_M$ is a metrizable, compact group with respect to $\Vert\cdot\Vert_2$.
Hence, if $\ep>0$ is fixed, there exists a finite set $1\in F\subset U_M$ such that 
\[
\inf\{\Vert u-v\Vert_2\colon u\in U_M, v\in F\}<\delta\coloneqq \frac{\ep}{6}.
\] 
Thus, if $\phi\in \CP_\tau(M)$ is subunital and such that
\[
\max_{v\in F}\Vert \phi(v)-v\Vert_2<\delta,
\]
then $\Vert \phi(u)-u\Vert_2\leq 3\delta$ for every $u\in U_M$, and as every $x\in (M)_1$ is of the form $x=\frac{1}{2}(v+v^*)+\frac{i}{2}(w+w^*)$ for suitable $v,w\in U_M$, we get that
\[
\sup_{x\in (M)_1}\Vert \phi(x)-x\Vert_2<\ep.
\]
This proves that $(M,\tau)$ satisfies condition (1) of Definition \ref{propT}. 
\hfill $\square$
\end{rem}

Finally, we give a characterization of II$_1$ factors with separable predual which do not have property (T) in terms of (relative) property gamma; it extends the case of group algebras as is presented in \cite[Theorem 4.6]{Jol3}.

More precisely, let $N\subset M$ be a pair of II$_1$ factors such that $N'\cap M=\C$, namely $N$ is \textit{irreducible} in $M$; recall from \cite{Jol3} that $N$ has \textit{property gamma relative to} $M$ if, for any free ultrafilter $\omega$ on $\N$, one has $N'\cap N^\omega\subsetneq N'\cap M^\omega$. Then we proved in \cite[Theorem 4.6]{Jol3} that if $G$ is an icc, non Kazhdan's group, there is an action of $G$ on the hyperfinite II$_1$ factor $R$ so that $L(G)$ is irreducible in the crossed product $R\rtimes G$ and $L(G)$ has property gamma relative to $R\rtimes G$. 

\begin{prop}
Let $N$ be a $\mathrm{II}_1$ factor with separable predual which does not have property (T). Then either $N$ has Murray and von Neumann property gamma or $N$ is full and there exists a $\mathrm{II}_1$ factor $M\supset N$ such that $N$ is irreducible in $M$ and $N'\cap M^\omega$ is diffuse. In particular, $N$ has property gamma relative to $M$.
\end{prop} 
\textsc{Proof.} If $N$ does not have property gamma, it is full, hence $N'\cap N^\omega=\C$. By \cite[Theorem 3.1]{Tan}, there is a tracial von Neumann algebra $(M,\tau)$ which contains $N$, and such that $N'\cap M=\C$, but such that $N'\cap M^\omega\not=\C$. It follows that $M$ is a II$_1$ factor in which $N$ is irreducible, and by \cite[Theorem 3.5]{FGL}, $N'\cap M^\omega$ is automatically diffuse.
\hfill $\square$

\vspace{1cm}
\noindent
\begin{flushright}
     \begin{tabular}{l}
       Universit\'e de Neuch\^atel,\\
       Institut de Math\'emathiques,\\       
       Emile-Argand 11\\
       CH-2000 Neuch\^atel, Switzerland\\
       \small {pajolissaint@gmail.com}
     \end{tabular}
\end{flushright}

\end{document}